# QUOTIENT SPACE OF $\mathcal{LMC}$-COMPACTIFICATION AS A SPACE OF $z-$FILTERS

M. AKBARI TOOTKABONI

ABSTRACT. The left multiplicative continuous compactification of a semitopological semigroup is the universal semigroup compactification. In this paper an internal construction of a semigroup compactification of a semitopological semigroup is constructed as a space of filters. In [6], we described an external construction of a semigroup compactification of a semitopological semigroup.

## 1. Introduction

Stone-Čech compactifications derived from a discrete semigroup S can be considered as the spectrum of the algebra $\mathcal{B}(S)$, the set of bounded complex-valued functions on $S$, or as a collection of ultrafilters on S. What is certain and indisputable is the fact that filters play an important role in the study of Stone-Čech compactifications derived from a discrete semigroup. It seems that filters can play a role in the study of general semigroup compactifications too, (See [6] and [7]).

For $S$ a semitopological semigroup, a continuous bounded function on $S$ is said to be in $\mathcal{LMC}(S)$ if its of right translates is relatively compact in $\mathcal{CB}(S)$ for the topology of pointwise convergence on $S$. In this paper we will begin with an elementary construction of semigroup compactification and we present some background about $\mathcal{LMC}(S)$. Also section 2 consists of an introduction to $z-$filters and an elementary external construction of semigroup compactification of a semitopological semigroup as a space of $z-$filters. In Glazer's proof of the finite sum theorem, he introduced the notion of the sum of two ultrafilter $p$ and $q$ on $\omega$, agreeing that $A \in p + q$ if and only if $\{x : -x + A \in q\} \in p$, (where by $-x + A$ is meant $\{y : x + y \in A\}$). In section 2, we introduce the notion of the binary operation of two $z-$filters.

Section three will be about some theorems from [2] about filters on discrete semigroup of S which are extended to a state that S be semeitopological semigroup. The scaffold of this and the following section focus on the coincidence on [2]. Therefore, all techniques and methods in [2] are established. In this section we investigate semigroup compactifications of $S$ as spaces of $z-$filters.

---

## 2. **Preliminaries**

Let $S$ be a semitopological semigroup (i.e. for each $s \in S$, $\lambda_s : S \to S$ and $r_s : S \to S$ are continuous, where for each $x \in S$, $\lambda_s(x) = sx$ and $r_s(x) = xs$) with a Hausdorff topology, and $\mathcal{CB}(S)$ denote the $C^*$-algebra of all bounded complex valued continuous functions on $S$. A semigroup compactification of $S$ is a pair $(\psi, X)$, where $X$ is a compact, Hausdorff, right topological semigroup (i.e. for all $x \in X$, $r_x$ is continuous) and $\psi : S \to X$ is continuous homomorphism with dense image such that, for all $s \in S$, the mapping $x \mapsto \psi(s)x : X \to X$ is continuous, (see Definition 3.1.1 in [1]). Let $\mathcal{F}$ be a $C^*$-subalgebra of $\mathcal{CB}(S)$ containing the constant functions, then the set of all multiplicative means of $\mathcal{F}$ (the spectrum of $\mathcal{F}$), which denote by $S^\mathcal{F}$, equipped with the Gelfand topology, is a compact Hausdorff topological space. A left translation invariant $C^*$-subalgebra $\mathcal{F}$ of $\mathcal{CB}(S)$ (i.e., $L_s f = f \circ \lambda_s \in \mathcal{F}$ for all $s \in S$ and $f \in \mathcal{F}$), containing the constant functions, is called $m$-admissible if the function $s \mapsto (T_\mu f(s)) = \mu(L_s f)$ is in $\mathcal{F}$ for all $f \in \mathcal{F}$ and $\mu \in S^\mathcal{F}$. If so, $S^\mathcal{F}$ under the multiplication $\mu\nu = \mu \circ T_\nu$ ($\mu, \nu \in S^\mathcal{F}$), furnished with the Gelfand topology, makes $(\varepsilon, S^\mathcal{F})$ a semigroup compactification (called the $\mathcal{F}$-compactification) of $S$, where $\varepsilon : S \to S^\mathcal{F}$ is the evaluation mapping. Also $\varepsilon^* : \mathcal{C}(S^\mathcal{F}) \to \mathcal{F}$ is an isometric isomorphism and $\widehat{f} = (\varepsilon^*)^{-1}(f) \in \mathcal{C}(S^\mathcal{F})$ for $f \in \mathcal{F}$ is given by $\widehat{f}(\mu) = \mu(f)$ for all $\mu \in S^\mathcal{F}$, (For more detail see section 2 in [1]). Now we present some famous m-admissible subalgebra of $\mathcal{CB}(S)$.

A function $f \in \mathcal{CB}(S)$ is left multiplicative continuous if and only if $\mathbf{T}_\mu f \in \mathcal{CB}(S)$ for all $\mu \in \beta S = S^{\mathcal{CB}(S)}$. We define

$$\mathcal{LMC}(S) = \bigcap \{\mathbf{T}_\mu^{-1}(\mathcal{CB}(S)) : \mu \in \beta S\}.$$

$\mathcal{LMC}(S)$ is the largest $m$- admissible subalgebra of $\mathcal{CB}(S)$. $(\varepsilon, S^{\mathcal{LMC}})$ is the universal compactification of $S$. ( Definition 4.5.1 and Theorem 4.5.2 in [1])

Now we quote some prerequisite material from [6] for the description of $(S^\mathcal{F}, \varepsilon)$ in terms of filters. For $f \in \mathcal{F}$, $Z(f) = f^{-1}(\{0\})$ is called zero set for all $f \in \mathcal{F}$ and we denote the collection of all zero sets by $Z(\mathcal{F})$. For an extensive account of ultrafilters, the readers may refer to [3], [4] and [5].

**Definition 2.1.** $\mathcal{A} \subseteq Z(\mathcal{F})$ is called a $z-$filter on $\mathcal{F}$ if
  (i) $\emptyset \notin \mathcal{A}$ and $S \in \mathcal{A}$,
  (ii) if $A, B \in \mathcal{A}$, then $A \bigcap B \in \mathcal{A}$,
  (iii) if $A \in \mathcal{A}$, $B \in Z(\mathcal{F})$ and $A \subseteq B$ then $B \in \mathcal{A}$.

A $z-$filter is said to be an $z-$ultrafilter if it is not contained properly in any other $z-$filters. We denote $\mathcal{F}S = \{p \subseteq Z(\mathcal{F}) : p \text{ is an } z-\text{ultrafilter on } \mathcal{F}\}$. It is obvious that $\widehat{x} = \{Z(f) : f \in \mathcal{F}, \ f(x) = 0\}$ is an $z-$ultrafilter on $\mathcal{F}$ and so $\widehat{x} \in \mathcal{F}S$. If $p, q \in \mathcal{F}S$, then obvious that
1) If $B \in Z(\mathcal{F})$ and for all $A \in p$, $A \cap B \neq \emptyset$ then $B \in p$,
2) If $A, B \in Z(\mathcal{F})$ such that $A \cup B \in p$ then $A \in p$ or $B \in p$ and

3) Let $p \neq q$ then there exist $A \in p$ and $B \in q$ such that $A \cap B = \emptyset$.( See Lemma 2.3 in [6]).

When $\mathcal{F} = \mathcal{CB}(S)$ then $\mathcal{F}S = \beta S = S^{\mathcal{CB}(S)}$ is the Stone-Čech compactification of $S$. Therefore the spectrum of $\mathcal{CB}(S)$ is equal with the collection of all $z$-ultrafilter on $S$. Now we describe semigroup compactification as the collection of $z-$filters subset of an $m-$ admissible subalgebra $\mathcal{F}$.

**Lemma 2.2.** *(i) Let $p \in \mathcal{F}S$, then there is a $\mu \in S^{\mathcal{F}}$ such that $\bigcap_{A \in p} \overline{\varepsilon(A)} = \{\mu\}$,*

*(ii) Let $p \in \mathcal{F}S$ and $\bigcap_{A \in p} \overline{\varepsilon(A)} = \{\mu\}$ for some $\mu \in S^{\mathcal{F}}$. If $A \in Z(\mathcal{F})$, and there exists a neighborhood $U$ of $\mu$ in $S^{\mathcal{F}}$ such that $U \subseteq \overline{\varepsilon(A)}$, then $A \in p$,*

*(iii) For each $\mu$ in $S^{\mathcal{F}}$ there exists a $p \in \mathcal{F}S$ such that $\bigcap_{A \in p} \overline{\varepsilon(A)} = \{\mu\}$.*

**Proof :** See Lemma 2.6 and 2.7 in [6]. $\square$

Unlike in the discrete setting, there is no simple correspondence between $S^{\mathcal{F}}$ and $\mathcal{F}S$. $\mathcal{F}S$ is equipped with a topology whose base is $\{(\widehat{A})^c : A \in Z(\mathcal{F})\}$, where $\widehat{A} = \{p \in \mathcal{F}S : A \in p\}$ is a compact space which is not Hausdorff in general. We define the relation $\sim$ on $\mathcal{F}S$ such that $p \sim q$ if
$$\bigcap_{A \in p} \overline{\varepsilon(A)} = \bigcap_{B \in q} \overline{\varepsilon(B)}.$$

It is obvious that $\sim$ is an equivalence relation on $\mathcal{F}S$. Let $[p]$ be the equivalence class of $p \in \mathcal{F}S$, and let $\frac{\mathcal{F}S}{\sim}$ be the corresponding quotient space with the quotient map $\pi : \mathcal{F}S \to \frac{\mathcal{F}S}{\sim}$. For every $p \in \mathcal{F}S$ define $\widetilde{p}$ by $\widetilde{p} = \bigcap [p]$, put $\widetilde{A} = \{\widetilde{p} : A \in p\}$ for $A \in Z(\mathcal{F})$ and $\mathcal{R} = \{\widetilde{p} : p \in \mathcal{F}S\}$.

**Lemma 2.3.** *The following statements hold:*

*(i) $\{(\widetilde{A})^c : A \in Z(\mathcal{F})\}$ is a basis for a topology on $\mathcal{R}$,*

*(ii) $\mathcal{R}$ is Hausdorff and compact,*

*(iii) The mapping $\varphi : S^{\mathcal{F}} \to \mathcal{R}$ defined by $\varphi(\mu) = \widetilde{p}$, in which $\bigcap_{A \in p} \overline{\varepsilon(A)} = \{\mu\}$, is a homeomorphism.*

**Proof :** See Lemma 2.11 and Lemma 2.12 in [6]. $\square$

Notice that by Lemma 2.6 , we have $\mathcal{R} = \{\mathcal{A}^{\mu} : \mu \in S^{\mathcal{F}}\}$, where $\mathcal{A}^{\mu} = \widetilde{p}$ for $\mu \in S^{\mathcal{F}}$ and $\bigcap_{A \in p} \overline{\varepsilon(A)} = \{\mu\}$, and also $\widehat{x} = \mathcal{A}^{\varepsilon(x)}$ , for each $x \in S$.
For all $x, y \in S$, we define
$$\mathcal{A}^{\varepsilon(x)} * \mathcal{A}^{\epsilon(y)} = \{Z(f) \in Z(\mathcal{F}) : Z(\mathbf{T}_{\varepsilon(y)}f) \in \mathcal{A}^{\varepsilon(x)}\}.$$

It is obvious that $\mathcal{A}^{\varepsilon(x)} * \mathcal{A}^{\varepsilon(y)} = \mathcal{A}^{\varepsilon(xy)}$, for each $x, y \in S$,( see Lemma 2.14 in [6]).

**Definition 2.4.** Let $\{x_{\alpha}\}$ and $\{y_{\beta}\}$ be two nets in $S$, such that $lim_{\alpha}\varepsilon(x_{\alpha}) = \mu$ and $lim_{\beta}\varepsilon(y_{\beta}) = \nu$, for $\mu, \nu \in S^{\mathcal{F}}$. We define
$$\mathcal{A}^{\mu} * \mathcal{A}^{\nu} = lim_{\alpha}(lim_{\beta}(\mathcal{A}^{\varepsilon(x_{\alpha})} * \mathcal{A}^{\varepsilon(y_{\beta})})).$$

**Theorem 2.5.** *The following statements hold:*

*(i) Definition 2.7 is well defined,*

*(ii) $(\mathcal{R}, e)$ is a compact right topological semigroup of $S$, in which $e : S \to \mathcal{R}$ is defined by $e[x] = \widehat{x}$,*

*(iii) The mapping $\varphi : S^{\mathcal{F}} \to \mathcal{R}$, defined by $\varphi(\mu) = \widetilde{p}$, where $\bigcap_{A \in p} \overline{\varepsilon(A)} = \{\mu\}$, is an isomorphism.*

**Proof :** See Theorem 2.16 in [6]. □

The operation "." on S extends uniquely to $(\mathcal{R}, *)$. Thus $(\mathcal{R}, e)$ is a semigroup compactification of $(S, .)$, where the evaluation map $e : S \to \mathcal{R}$ is given $e[x] = \widehat{x}$, for each $x \in S$. Also, $e[S]$ is a subset of the topological center of $\mathcal{R}$ and $cl_{\mathcal{R}}(e[S]) = \mathcal{R}$. For more details see [6]. So in this paper we denote $S^{\mathcal{F}} = \{\mathcal{A}^{\mu} : \mu \in S^{\mathcal{F}}\}$, where $\bigcap_{A \in \mathcal{A}^{\mu}} \overline{\varepsilon(A)} = \{\mu\}$.

**Definition 2.6.** A $z-$filter $\mathcal{A}$ on $\mathcal{F}$ is called a pure $z-$filter if for some $p \in \mathcal{F}S$, $\mathcal{A} \subseteq p$ implies that $\mathcal{A} \subseteq \widetilde{p}$. The collection of all pure $z-$filters are denoted by $\mathcal{P}(\mathcal{F})$.

It is obvious that $\widetilde{p} \in S^{\mathcal{F}}$ is a maximal member of $\mathcal{P}(\mathcal{F})$.

**Definition 2.7.** For a $z-$filter $\mathcal{A} \subseteq Z(\mathcal{F})$, we define

(i) $\overline{\mathcal{A}} = \{\widetilde{p} \in S^{\mathcal{F}} : \text{There exists } p \in \mathcal{F}S \text{ such that } \mathcal{A} \subseteq p\}$,

(ii) $\mathcal{A}^{\circ} = \{A \in \mathcal{A} : \overline{\mathcal{A}} \subseteq (\overline{\varepsilon(A)})^{\circ}\}$.

**Lemma 2.8.** *Let $\mathcal{A}$ and $\mathcal{B}$ be $z-$filters. Then the following statements hold.*

*(1) $\overline{\mathcal{A}} = \bigcap_{A \in \mathcal{A}^{\circ}} (\overline{\varepsilon(A)})^{\circ} = \bigcap_{A \in \mathcal{A}^{\circ}} \overline{\varepsilon(A)} = \bigcap_{A \in \mathcal{A}} \overline{\varepsilon(A)}$.*

*(2) Let $\mathcal{A} \subseteq Z(\mathcal{F})$ be a $z-$filter, then $\overline{\mathcal{A}}$ is a closed subset of $S^{\mathcal{F}}$.*

*(3) $\mathcal{A} \subseteq \mathcal{B}$ then $\overline{\mathcal{B}} \subseteq \overline{\mathcal{A}}$.*

*(4) Let $\widetilde{p} \in S^{\mathcal{F}}$. If $A \in Z(\mathcal{F})$ and $\widetilde{p} \in (\overline{\varepsilon(A)})^{\circ}$, then $A \in \widetilde{p}$.*

*(5) If $\mathcal{J}$ is a closed subset of $S^{\mathcal{F}}$ then $\mathcal{A} = \bigcap \mathcal{J}$ is a pure $z-$filter and $\overline{\mathcal{A}} = \mathcal{J}$.*

*(6) Let $\mathcal{A}$ be a $z-$filter, then $\bigcap \overline{\mathcal{A}}$ is a pure $z-$filter and $\bigcap \overline{\mathcal{A}} \subseteq \mathcal{A}$. In addition, if $\mathcal{A}$ is a pure $z-$filter then $\mathcal{A} = \bigcap \overline{\mathcal{A}}$.*

*(7) Let $\mathcal{A}$ and $\mathcal{B}$ be pure $z-$filters then $\mathcal{A} \subseteq \mathcal{B}$ if and only if $\overline{\mathcal{B}} \subseteq \overline{\mathcal{A}}$, also $\mathcal{A} = \mathcal{B}$ if and only if $\overline{\mathcal{A}} = \overline{\mathcal{B}}$.*

*(8) Let $A, B \in Z(\mathcal{F})$. Then $(\overline{\varepsilon(A)})^{\circ} \cap (\overline{\varepsilon(B)})^{\circ} = (\overline{\varepsilon(A \cap B)})^{\circ}$.*

*(9) $\mathcal{A}^{\circ}$ is a $z-$filter.*

*(10) Let $\mathcal{A} \in \mathcal{F}S$, $A \in Z(\mathcal{F})$ and $A \notin \mathcal{A}$. If $F \in Z(\mathcal{F})$ and $A^c \subseteq F$, then $F \in \mathcal{A}$.*

*(11) Let $\widetilde{p} \in S^{\mathcal{F}}$, $A \in Z(\mathcal{F})$ and $A \notin \widetilde{p}$. If $F \in Z(\mathcal{F})$ and $A^c \subseteq F$, then $\widetilde{p} \in \overline{\varepsilon(F)}$.*

*(12) Let $p \in \mathcal{F}S$, $A \in Z(\mathcal{F})$ and $\widetilde{p} \notin \overline{\varepsilon(A)}$. If $F \in Z(\mathcal{F})$ and $A^c \subseteq F$, then $F \in \widetilde{p}$.*

**Proof :** (1), (2) and (3) are straightforward.

(4) By Lemma 2.5(ii), it is obvious.

(5) It is obvious that $\mathcal{A} = \bigcap \mathcal{J}$ is a pure $z-$filter and $\mathcal{J} \subseteq \overline{\mathcal{A}}$. Let $f \in \mathcal{F}$ such that $\mathcal{J} \subseteq (\overline{\varepsilon(Z(f))})^{\circ}$. Since for each $\widetilde{p} \in \mathcal{J} \subseteq (\overline{\varepsilon(Z(f))})^{\circ}$ we have $Z(f) \in \widetilde{p}$ (Lemma

3.2), so $Z(f) \in \mathcal{A}$. Hence by $(i)$
$$\overline{\mathcal{A}} = \bigcap_{A \in \mathcal{A}} \overline{A} \subseteq \bigcap_{\mathcal{J} \subseteq (\overline{\varepsilon(A)})^\circ} \overline{A} \subseteq \mathcal{J},$$
and therefore $\overline{\mathcal{A}} = \mathcal{J}$.

(6) It is obvious that $\bigcap \overline{\mathcal{A}}$ is a pure $z-$filter and $\bigcap \overline{\mathcal{A}} \subseteq \mathcal{A}$. Now let $\mathcal{A}$ be a pure $z-$filter, then for every $\widetilde{p} \in \overline{\mathcal{A}}$ we have $\mathcal{A} \subseteq \widetilde{p}$, and this implies that $\mathcal{A} \subseteq \bigcap_{\widetilde{p} \in \overline{\mathcal{A}}} \widetilde{p} = \bigcap \overline{\mathcal{A}}$.

(7) is straightforward.

(8) If $A, B \in Z(\mathcal{F})$ then $\varepsilon(A \cap B) = \varepsilon(A) \cap \varepsilon(B)$. $(\overline{\varepsilon(A \cap B)})^\circ \subseteq (\overline{\varepsilon(A)})^\circ \cap (\overline{\varepsilon(B)})^\circ$ follows from basic topology.

For the converse, let $x \in (\overline{\varepsilon(A)})^\circ \cap (\overline{\varepsilon(B)})^\circ$, then there exists $f \in \mathcal{F}$ such that $x \in (\overline{\varepsilon(Z(f))})^\circ \subseteq \overline{\varepsilon(Z(f))} \subseteq (\overline{\varepsilon(A)})^\circ \cap (\overline{\varepsilon(B)})^\circ$. Hence $Z(f) \subseteq \varepsilon(A \cap B)$ and so $x \in (\overline{\varepsilon(Z(f))})^\circ \subseteq \overline{\varepsilon(Z(f))} \subseteq \overline{\varepsilon(A \cap B)}$ and this implies that $x \in (\overline{\varepsilon(A \cap B)})^\circ$.

(9) is straightforward.

(10) Let $F \in Z(\mathcal{F})$ be such that $A^c \subseteq F$. Then $\mathcal{A} \cup \{F\}$ has the finite intersection property.( Let there exist $A_1, ..., A_n \in \mathcal{A}$ such that $(\cap_{i=1}^n A_i) \cap F = \emptyset$ then $\cap_{i=1}^n A_i \subseteq F^c \subseteq A$ and so $A \in \mathcal{A}$ is a contradiction.) Therefore there exists $\mathcal{B} \in \mathcal{FS}$ such that $\mathcal{A} \cup \{F\} \subseteq \mathcal{B}$. Since $\mathcal{A}$ is an $z-$ultrafilter so $\mathcal{A} = \mathcal{B}$, and this implies that $F \in \mathcal{A}$.

(11) Let $A \notin \widetilde{p}$, so there exists $q \in [p]$ such that $A \notin q$. Hence by (10), for every $F \in Z(\mathcal{F})$ if $A^c \subseteq F$ then $F \in q$. Therefore $\widetilde{p} \in \overline{\varepsilon(F)}$.

(12) Let $\widetilde{p} \notin \overline{\varepsilon(A)}$ then $[p] \cap \overline{\varepsilon(A)} = \emptyset$. If $q \in [p]$ then $A \notin q$. Hence by (10), for each $F \in Z(\mathcal{F})$ and for each $q \in [p]$ if $A^c \subseteq F$ then $F \in q$. Therefore $F \in \widetilde{p} = \bigcap_{q \in [p]} q$. $\square$

**Definition 2.9.** Let $A$ be a nonempty closed subset of $S^\mathcal{F}$. We call $\bigcap A$ the pure $z-$filter generated by $A$.

**Definition 2.10.** Let $\mathcal{A}$ and $\mathcal{B}$ be two $z-$filters on $\mathcal{F}$ and $A \in Z(\mathcal{F})$. We say $A \in \mathcal{A} + \mathcal{B}$ if and only if for each $F \in Z(\mathcal{F})$, $\Omega_\mathcal{B}(A) \subset F$ implies $F \in \mathcal{A}$, where $\Omega_\mathcal{B}(A) = \{x \in S : \lambda_x^{-1}(A) \in \mathcal{A}\}$.

**Lemma 2.11.** Let $\mathcal{A}$ and $\mathcal{B}$ be two $z-$filters and $A, B \in Z(\mathcal{F})$. Then

(1) $\Omega_\mathcal{A}(A \cap B) = \Omega_\mathcal{A}(A) \cap \Omega_\mathcal{A}(B)$,

(2) If $\mathcal{A} \subseteq \mathcal{B}$ then $\Omega_\mathcal{A}(A) \subseteq \Omega_\mathcal{B}(A)$,

(3) If $A \subseteq B$ then $\Omega_\mathcal{A}(A) \subseteq \Omega_\mathcal{A}(B)$,

(4) For each $x \in S$, $\lambda_x^{-1}(\Omega_\mathcal{A}(A)) = \Omega_\mathcal{A}(\lambda_x^{-1}(A))$.

(5) $\Omega_{\mathcal{B}+\mathcal{C}}(A) = \Omega_\mathcal{B}(\Omega_\mathcal{C}(A))$

(6) Let $\widetilde{p}, \widetilde{q} \in S^\mathcal{F}$ and $B \in Z(\mathcal{F})$. If $\widetilde{q} * \widetilde{p} \in (\overline{\varepsilon(B)})^\circ$ then
$$\varepsilon^{-1}(r_{\widetilde{p}}^{-1}((\overline{\varepsilon(B)})^\circ)) \subseteq \{x \in S : \widetilde{p} \in (\overline{\varepsilon(\lambda_x^{-1}(B))})^\circ\} \subseteq \Omega_{\widetilde{p}}(B).$$

**Proof :** Obvious. $\square$

**Lemma 2.12.** *(i) Let $\widetilde{p} \in S^{\mathcal{LMC}}$ and $A \in Z(\mathcal{LMC}(S))$. Then $\Omega_{\widetilde{p}}(A)$ is a closed subset of $S$.*
*(ii) Let $\mathcal{B}$ be a $z-$filter and $A \in Z(\mathcal{LMC}(S))$. Then $\Omega_{\mathcal{B}}(A)$ is a closed subset of $S$.*

**Proof :** Let $\widetilde{p} \in S^{\mathcal{LMC}}$ and $A \in Z(\mathcal{LMC}(S))$, then

$$
\begin{aligned}
\Omega_{\widetilde{p}}(A) &= \{x \in S : \lambda_x^{-1}(A) \in \widetilde{p}\} \\
\text{for some } f \in \mathcal{LMC}(S) \quad &= \{x \in S : Z(L_x f) \in \widetilde{p}\} \\
&= \bigcap_{p \in \beta S, \widetilde{p} \subseteq p} \{x \in S : Z(L_x f) \in p\} \\
\text{for a zero set } B \text{ and } p \in \beta S, B \in p \text{ iff } p \in \overline{B} \quad &= \bigcap_{p \in \beta S, \widetilde{p} \subseteq p} \{x \in S : p \in \overline{Z(L_x f)}\} \\
&= \bigcap_{\mu \in \beta S, \widetilde{p} \subseteq \mathcal{A}^\mu} \{x \in S : \mathbf{T}_\mu(f)(x) = 0, \} \\
&= \bigcap_{\mu \in \beta S, \widetilde{p} \subseteq \mathcal{A}^\mu} Z(\mathbf{T}_\mu(f)).
\end{aligned}
$$

Since $f \in \mathcal{LMC}(S)$, hence $\mathbf{T}_\mu(f) : S \to \mathbb{C}$ is continuous and $Z(\mathbf{T}_\mu(f))$ is a closed subset of $S$. Hence $\Omega_{\widetilde{p}}(A)$ is a closed subset of $S$.

$(ii)$ It is obvious. □

**Lemma 2.13.** *Let $\mathcal{A}$ and $\mathcal{B}$ be $z-$filters on $\mathcal{LMC}(S)$ then $\mathcal{A} + \mathcal{B}$ is a $z-$filter.*

**Proof :** It is obvious $\emptyset \notin \mathcal{A} + \mathcal{B}$. Let $A, B \in \mathcal{A} + \mathcal{B}$ and $F \in Z(\mathcal{LMC}(S))$ such that $\Omega_{\mathcal{B}}(A \cap B) \subseteq F$ then there exist $D, C \in Z(\mathcal{LMC}(S))$ such that $\Omega_{\mathcal{B}}(A) \subset C$, $\Omega_{\mathcal{B}}(B) \subset D$ and $C \cap D \subseteq F$. This implies $F \in \mathcal{A}$. Now let $A \in \mathcal{A} + \mathcal{B}$ and $A \subseteq B$ for some $B \in Z(\mathcal{LMC}(S))$, then by Lemma 2.14(3) implies $B \in \mathcal{A} + \mathcal{B}$. □

**Definition 2.14.** Let $\mathcal{A}$ and $\mathcal{B}$ be $z-$filters, we define $\mathcal{A} \odot \mathcal{B} = \bigcap \overline{\mathcal{A} + \mathcal{B}}$.

By Definition 2.12, $\mathcal{A} \odot \mathcal{B}$ is a pure $z-$filter generated by $\overline{\mathcal{A} + \mathcal{B}}$.

**Lemma 2.15.** *Let $\mathcal{A}$ and $\mathcal{B}$ be $z-$filters and $\mu, \nu \in S^{\mathcal{LMC}}$, then*
*(i) $\mathcal{A} \odot \mathcal{A}^\mu = \bigcap_{\mathcal{A}^\nu \in \overline{\mathcal{A}}} \mathcal{A}^{\nu\mu}$,*
*(ii) $\mathcal{A}^\mu \odot \mathcal{A}^\nu = \mathcal{A}^{\mu\nu} = \mathcal{A}^\mu * \mathcal{A}^\nu$,*
*(iii) $\mathcal{A} \odot \mathcal{B} \subseteq \bigcap_{\mathcal{A}^\nu \in \overline{\mathcal{B}}, \mathcal{A}^\mu \in \overline{\mathcal{A}}} \mathcal{A}^{\mu\nu}$,*
*(iv) $\mathcal{A}^\mu \odot \mathcal{A} \subseteq \bigcap_{\mathcal{A}^\nu \in \overline{\mathcal{A}}} \mathcal{A}^{\mu\nu}$.*
*(v) $\overline{\mathcal{A}} * \overline{\mathcal{B}} \subseteq cl_{S^{\mathcal{F}}}(\overline{\mathcal{A}} * \overline{\mathcal{B}}) \subseteq \overline{\mathcal{A} + \mathcal{B}} = \overline{\mathcal{A} \odot \mathcal{B}}.$*

**Proof :** Obvious. □

In this paper, we replace $\overline{\varepsilon(A)}$ with $\overline{A}$ for simplicity.

## 3. Compactification of $S$ as a space of pure $z-$filters

Recall from Theorem 4.2.1 in [1] that if $S$ is a semitopological semigroup, $(S^{\mathcal{LMC}}, \varepsilon)$ is the maximal right topological semigroup extention of $S$. Thus given any semigroup compactification $(T, f)$ of $S$, we have $\widehat{f}$ is a continuous homomorphism of $S^{\mathcal{LMC}}$ onto $T$. Consequently, $(T,.)$ is a topological quotient of $S^{\mathcal{LMC}}$. ( Let $R(f) = \{(\mu,\nu) : \mu,\nu \in S^{\mathcal{LMC}}$ and $\widehat{f}(\mu) = \widehat{f}(\nu)\}$. Then $T \approx S^{\mathcal{LMC}}/R(f)$ via the homomorphism $g([\mu]) = f(\mu)$. Furthermore, defining $\oplus$ on $S^{\mathcal{LMC}}/R(f)$ by $[\mu\nu] = [\mu] \oplus [\nu]$ we have $g$ is an isomorphism of $S^{\mathcal{LMC}}/R(f)$ onto $(T,.)$.)

Given an equivalence relation $R$ on $S^{\mathcal{LMC}}$ such that $S^{\mathcal{LMC}}/R(f)$ is Hausdorff, we have the $R$-equivalence class are closed subsets of $S^{\mathcal{LMC}}$. But the closed subsets of $S^{\mathcal{LMC}}$ correspond exactly to the pure $z$−filters.( See Lemma 2.11(5).)

We thus have that any semigroup compactification of $S$ corresponds to a set of pure $z$−filters. In this section we characterize those sets of pure $z$−filters which yield semigroup compactification of $S$ and describe the topology and the operation in terms of the $z$−filters.

**Definition 3.1.** Let $\Gamma$ be a subset of pure $z$−filters on $\mathcal{LMC}(S)$. We say $f : \Gamma \to \bigcup \Gamma$ is a topological choice function if for each $\mathcal{L} \in \Gamma$, $\overline{\mathcal{L}} \subseteq (\overline{f(\mathcal{L})})^\circ$.

**Theorem 3.2.** *Let $\Gamma$ be a set of pure $z$−filters on $\mathcal{LMC}(S)$. Statements $(a)$ and $(a')$ are equivalent, statements $(b)$ and $(b')$ are equivalent and statement $(c)$ is equivalent to the conjunction of statements $(a)$ and $(b)$.*

*(a) Given any topological choice function $f$ for $\Gamma$, there is a finite subfamily $\mathcal{F}$ of $\Gamma$ such that $\overline{S} = S^{\mathcal{LMC}} = \bigcup_{\mathcal{L} \in \mathcal{F}} (\overline{f(\mathcal{L})})^\circ$.*

*(a') For each $\widetilde{p} \in S^{\mathcal{LMC}}$ there is some $\mathcal{L} \in \Gamma$ such that $\mathcal{L} \subseteq \widetilde{p}$.*

*(b) Given distinct $\mathcal{L}$ and $\mathcal{K}$ in $\Gamma$, there exists $B \in \mathcal{K}^\circ$ such that for every $A \in Z(\mathcal{LMC}(S))$ if $\overline{S} - (\overline{B})^\circ \subseteq (\overline{A})^\circ$ then $A \in \mathcal{L}$.*

*(b') For each $\widetilde{p} \in S^{\mathcal{LMC}}$ there is at most one $\mathcal{L} \in \Gamma$ such that $\mathcal{L} \subseteq \widetilde{p}$.*

*(c) There is an equivalence relation $R$ on $S^{\mathcal{LMC}}$ such that each equivalence class is closed in $S^{\mathcal{LMC}}$ and $\Gamma = \{\, \bigcap [\widetilde{p}]_R : \widetilde{p} \in S^{\mathcal{LMC}}\}$.*

**Proof :** To see that $(a')$ implies $(a)$ let $f$ be a topological choice function for $\Gamma$. Suppose that the conclusion of $(a)$ fails. Then $\Gamma = \{\overline{S} - (\overline{f(\mathcal{L})})^\circ : \mathcal{L} \in \Gamma\}$ has the finite intersection property, therefore there exists $\widetilde{p} \in \bigcap_{T \in \Gamma} T$. Pick $\mathcal{L} \in \Gamma$ such that $\mathcal{L} \subseteq \widetilde{p}$. Since $f$ is a topological choice function $\widetilde{p} \in \overline{\mathcal{L}} \subseteq (\overline{f(\mathcal{L})})^\circ$ and so $\widetilde{p} \in (\overline{f(\mathcal{L})})^\circ \bigcap (\bigcap_{\mathcal{K} \in \Gamma} \overline{S} - (\overline{f(\mathcal{K})})^\circ) \subseteq (\overline{f(\mathcal{L})})^\circ \bigcap (\overline{S} - (\overline{f(\mathcal{L})})^\circ) = \emptyset$ is a contradiction.

To see that $(a)$ implies $(a')$ let $\widetilde{p} \in S^{\mathcal{LMC}}$ and suppose that the conclusion of $(a')$ fails. For each $\mathcal{L} \in \Gamma$, pick $f(\mathcal{L}) \in \mathcal{L} - \widetilde{p}$ such that $\overline{\mathcal{L}} \subseteq (\overline{f(\mathcal{L})})^\circ$. Pick a finite subfamily $\mathcal{F}$ of $\Gamma$ such that $S^{\mathcal{LMC}} = \bigcup_{\mathcal{L} \in \mathcal{F}} (\overline{f(\mathcal{L})})^\circ$. Then for some $\mathcal{L} \in \mathcal{F}$, $f(\mathcal{L}) \in \widetilde{p}$, a contradiction.

To see that $(b')$ implies $(b)$ let $\mathcal{L}$ and $\mathcal{K}$ be distinct members of $\Gamma$ and suppose that for each $B \in \mathcal{K}^\circ$ there exists $A \in Z(\mathcal{LMC}(S))$ such that $\overline{S} - (\overline{B})^\circ \subseteq (\overline{A})^\circ$ and $A \notin \mathcal{L}$. Then for each $B \in \mathcal{K}^\circ$, $\mathcal{L} \cup \{B\}$ has the finite intersection property. (Let $\mathcal{F}$ be a finite subfamily of $\mathcal{L}$ and $B \bigcap (\cap \mathcal{F}) = \emptyset$ then $(\cap \mathcal{F}) \bigcap (\overline{B})^\circ \subseteq (\cap \mathcal{F}) \bigcap (\overline{B}) = \emptyset$. Hence $\cap \mathcal{F} \subseteq \overline{S} - (\overline{B})^\circ \subseteq (\overline{A})^\circ$ implies $\overline{\mathcal{L}} \subseteq \overline{\cap \mathcal{F}} \subseteq [(\overline{B})^\circ]^c \subseteq (\overline{A})^\circ$ and so $A \in \mathcal{L}$, and this is a contradiction.) Since $\mathcal{K}$ is closed under finite intersection , $\mathcal{L} \bigcup \mathcal{K}^\circ$ has the finite intersection property. But then one can pick $\widetilde{p} \in S^{\mathcal{LMC}}$ such that $\mathcal{L} \bigcup \mathcal{K}^\circ \subseteq \widetilde{p}$ and this is contradiction $(b')$.

To see that $(b)$ implies $(b')$, let $\widetilde{p} \in S^{\mathcal{LMC}}$ and suppose we have distinct $\mathcal{L}$ and $\mathcal{K}$ in $\Gamma$ such that $\mathcal{L} \bigcup \mathcal{K} \subseteq \widetilde{p}$. Pick $B \in \mathcal{K}^\circ$. Let $A \in Z(\mathcal{LMC}(S))$ and $\overline{S} - (\overline{B})^\circ \subseteq (\overline{A})^\circ$ then $A \in \mathcal{L}$, and so $\widetilde{p} \in (\overline{B})^\circ \bigcap (\bigcap_{\overline{S}-(\overline{B})^\circ \subseteq (\overline{A})^\circ} \overline{A}) = (\overline{B})^\circ \bigcap (\overline{S} - (\overline{B})^\circ) = \emptyset$. Hence we have a contradiction.

That $(c)$ implies $(a')$ is trivial since $\bigcap [\widetilde{p}] \subseteq \widetilde{p}$.

To see that $(c)$ implies $(b')$ suppose $\widetilde{q} \in S^{\mathcal{LMC}}$ such that $[\widetilde{q}] \neq [\widetilde{p}]$ and $\bigcap [\widetilde{q}] \subseteq \widetilde{p}$. Since $[\widetilde{q}]$ is closed in $S^{\mathcal{LMC}}$ and $\widetilde{p} \notin [\widetilde{q}]$, pick $A \in (\widetilde{p})^\circ$ such that $\overline{A} \bigcap [\widetilde{q}] = \emptyset$. Let $B \in Z(\mathcal{LMC}(S))$ and $[\widetilde{q}] \subseteq \overline{S} - (\overline{A})^\circ \subseteq (\overline{B})^\circ$, then $B \in \bigcap [\widetilde{q}] \subseteq \widetilde{p}$ and so this is a contradiction.

Now assume that $(a')$ and $(b')$ hold. Define $R$ on $S^{\mathcal{LMC}}$ by $\widetilde{p}R\widetilde{q}$ if and only if there exists $\mathcal{L} \in \Gamma$ such that $\mathcal{L} \subseteq \widetilde{p} \bigcap \widetilde{q}$. Trivially R is symmetric. By $(a')$, R is reflexive on $S^{\mathcal{LMC}}$. To see that $R$ is transitive, let $\widetilde{p}, \widetilde{q}$ and $\widetilde{r}$ be in $S^{\mathcal{LMC}}$ and assume that $\widetilde{p}R\widetilde{q}$ and $\widetilde{q}R\widetilde{r}$. Pick $\mathcal{L}$ and $\mathcal{K}$ in $\Gamma$ such that $\mathcal{L} \subseteq \widetilde{p} \bigcap \widetilde{q}$ and $\mathcal{K} \subseteq \widetilde{q} \bigcap \widetilde{r}$. Since $\mathcal{L} \bigcup \mathcal{K} \subseteq \widetilde{q}$ we have $\mathcal{L} = \mathcal{K}$ by $(b')$.

To see that $\Gamma = \{ \bigcap [\widetilde{p}]_R : \widetilde{p} \in S^{\mathcal{LMC}} \}$, it suffices to show that if $\mathcal{L} \in \Gamma$ and $\widetilde{p} \in S^{\mathcal{LMC}}$ such that $\mathcal{L} \subseteq \widetilde{p}$, then $\mathcal{L} = \bigcap [\widetilde{p}]$. (For given any $\mathcal{L} \in \Gamma$ there is some $\widetilde{p} \in S^{\mathcal{LMC}}$ such that $\mathcal{L} \subseteq \widetilde{p}$ and given $\widetilde{p} \in S^{\mathcal{LMC}}$ we have by $(a')$ some $\mathcal{L} \in \Gamma$ such that $\mathcal{L} \subseteq \widetilde{p}$.) Let $\mathcal{L} \in \Gamma$ and $\widetilde{p} \in S^{\mathcal{LMC}}$ with $\mathcal{L} \subseteq \widetilde{p}$ be given. Given $\widetilde{q} \in [\widetilde{p}]$ there is some $\mathcal{K} \in \Gamma$ such that $\mathcal{K} \subseteq \widetilde{p} \bigcap \widetilde{q}$. But then, $\mathcal{K} = \mathcal{L}$ by $(b')$, so $\mathcal{L} \subseteq \widetilde{q}$. Thus $\mathcal{L} \subseteq \bigcap [\widetilde{p}]$.

To see that $\cap [\widetilde{p}] \subseteq \mathcal{L}$, let $A \in \cap [\widetilde{p}]$ such that $[\widetilde{p}] \subseteq (\overline{A})^\circ$ and $\overline{\mathcal{L}} - (\overline{A}) \neq \emptyset$. We choice $B \in Z(\mathcal{LMC}(S))$ such that $\overline{S} - (\overline{A})^\circ \subseteq (\overline{B})^\circ$ and $B \notin \mathcal{L}$. Then $\mathcal{L} \bigcup \{D\}$ has the finite intersection property for each $D \in Z(\mathcal{LMC}(S))$ such that $\overline{S} - (\overline{A})^\circ \subseteq (\overline{D})^\circ \subseteq \overline{D} \subseteq (\overline{B})^\circ$, (If $\mathcal{F}$ is a finite subfamily of $\mathcal{L}$ and $D \bigcap (\cap \mathcal{F}) = \emptyset$ then $(\cap \mathcal{F}) \bigcap (\overline{D})^\circ \subseteq (\cap \mathcal{F}) \bigcap (\overline{D}) = \emptyset$. Let $D \in Z(\mathcal{LMC}(S))$, if $\overline{S} - (\overline{A})^\circ \subseteq (\overline{D})^\circ \subseteq \overline{D} \subseteq (\overline{B})^\circ$, then $(\bigcap \mathcal{F}) \bigcap (\overline{D})^\circ = \emptyset$ and $(\cap \mathcal{F}) \bigcap (\overline{S} - (\overline{A})^\circ) = \emptyset$. Hence $\cap \mathcal{F} \subseteq (\overline{A})^\circ$ and so $\overline{\mathcal{L}} - \overline{A} = \emptyset$ and this is a contradiction). So pick $\widetilde{q} \in S^{\mathcal{LMC}}$ such that $\widetilde{q} \in \overline{\mathcal{L} \bigcup \{B\}}$. But then $\mathcal{L} \subseteq \widetilde{p} \bigcap \widetilde{q}$ so $\widetilde{q} \in [\widetilde{p}]$ and hence $\widetilde{q} \in (\overline{A})^\circ$. Also for each $E \in Z(\mathcal{LMC}(S))$, if $\overline{S} - (\overline{A})^\circ \subseteq (\overline{E})^\circ \subseteq (\overline{E}) \subseteq (\overline{B})^\circ$ then $E \in \widetilde{q}$. So $\widetilde{q} \in \overline{A} \bigcap (\overline{S} - (\overline{A})^\circ)$, and this is a contradiction.

Finally to show that each $R$- equivalence class is closed in $S^{\mathcal{LMC}}$, let $\widetilde{p} \in S^{\mathcal{LMC}}$ and pick $\mathcal{L} \in \Gamma$ such that $\mathcal{L} \subseteq \widetilde{p}$. Then $\overline{\mathcal{L}} = [\widetilde{p}]$ and, as we have observed $\overline{\mathcal{L}}$ is closed in $S^{\mathcal{LMC}}$. $\square$

**Definition 3.3.** Let $\Gamma \subseteq Z(\mathcal{LMC}(S))$ be a set of pure $z-$filters. We define

$$A^* = \{\mathcal{L} \in \Gamma : \overline{\mathcal{L}} \bigcap \overline{A} \neq \emptyset\},$$
$$A_* = \{\mathcal{L} \in \Gamma : \overline{\mathcal{L}} \subseteq \overline{A}\}$$

for each $A \in Z(\mathcal{LMC}(S))$.

Let $\Gamma$ be a collection of $z-$filters and let $A, B \in Z(\mathcal{LMC}(S))$, then we have

$$
\begin{aligned}
(i) & \quad (A \cap B)^* \subseteq A^* \cap B^*, \\
(ii) & \quad A^* \cup B^* = (A \cup B)^*, \\
(iii) & \quad (A \cap B)_* \subseteq A_* \cap B_*, \\
(iv) & \quad (A \cup B)_* = A_* \cup B_*.
\end{aligned}
$$

Therefore $\{(A^*)^c : A \in Z(\mathcal{LMC}(S))\}$ and $\{(A_*)^c : A \in Z(\mathcal{LMC}(S))\}$ are basis for the topology on $\Gamma$. It is obvious that $A_* \subseteq A^*$ and so $(A^*)^c \subseteq (A_*)^c$ for each $A \in Z(\mathcal{LMC}(S))$, and this is a motivation for the following definition.

**Definition 3.4.** Let $\Gamma \subseteq Z(\mathcal{F})$ be a set of pure $z-$filters. The $\tau^*$ quotient topology on $\Gamma$ is the topology generated by $\{(A^*)^c : A \in Z(\mathcal{LMC}(S))\}$ and the $\tau_*$ quotient topology on $\Gamma$ is the topology generated by $\{(A_*)^c : A \in Z(\mathcal{LMC}(S))\}$.

Let $R$ be an equivalence relation on $S^{\mathcal{LMC}}$ and let $\Gamma = \{\cap [\widetilde{p}]_R : \widetilde{p} \in S^{\mathcal{LMC}}\}$. Then for every $A \in Z(\mathcal{LMC}(S))$ we have $A^* = \{\cap [\widetilde{p}]_R : \widetilde{p} \in \overline{A}\}$ and $A_* = \{\cap [\widetilde{p}]_R : [\widetilde{p}] \subseteq \overline{A}\}$.

**Lemma 3.5.** *Let $R$ be an equivalence relation on $S^{\mathcal{LMC}}$ such that $S^{\mathcal{LMC}}/R$ is Hausdorff. Let $\Gamma = \{\cap [\widetilde{p}]_R : \widetilde{p} \in S^{\mathcal{LMC}}\}$. Then with the $\tau^*$ quotient topology on $\Gamma$, $\Gamma \approx S^{\mathcal{LMC}}/R$.*

**Proof :** Define $\varphi : S^{\mathcal{LMC}}/R \to \Gamma$ by $\varphi([\widetilde{p}]) = \cap [\widetilde{p}]$. It is obvious that $\varphi$ is well defined, onto and one to one. Let $\pi : S^{\mathcal{LMC}} \to S^{\mathcal{LMC}}/R$ is the quotient map and let $\gamma : S^{\mathcal{LMC}} \to \Gamma$ is defined by $\gamma(\widetilde{p}) = \cap [\widetilde{p}]$. It is obvious that $\gamma$ is well defined and onto. Also $\varphi \circ \pi = \gamma$, i.e. the following diagram

$$
\begin{array}{ccc}
S^{\mathcal{LMC}} & \xrightarrow{\pi} & S^{\mathcal{LMC}}/R \\
{\scriptstyle \gamma} \searrow & & \swarrow {\scriptstyle \varphi} \\
& \Gamma &
\end{array}
$$

commutes.

$\gamma$ is continuous, because for each $A \in Z(\mathcal{LMC}(S))$, we have

$$
\begin{aligned}
\gamma^{-1}(A^*) &= \{\widetilde{p} \in S^{\mathcal{LMC}} : \gamma(\widetilde{p}) \in A^*\} \\
&= \{\widetilde{p} \in S^{\mathcal{LMC}} : \varphi \circ \pi(\widetilde{p}) \in A^*\} \\
&= \{\widetilde{p} \in S^{\mathcal{LMC}} : \cap [\widetilde{p}] \in A^*\} \\
&= \{\widetilde{p} \in S^{\mathcal{LMC}} : \overline{A} \cap [\widetilde{p}] \neq \emptyset\} \\
&= \bigcup \{\pi(\widetilde{p}) : \widetilde{p} \in \overline{A}\} \\
&= \pi^{-1}(\pi(\{\widetilde{p} \in S^{\mathcal{LMC}} : \widetilde{p} \in \overline{A}\})) \\
&= \pi^{-1}(\pi(\overline{A})).
\end{aligned}
$$

Since $\overline{A}$ is compact and $\pi$ is continuous, $\pi(\overline{A})$ is compact, and so it is closed. Hence $\pi^{-1}(\pi(\overline{A})) = \gamma^{-1}(A^*)$ is closed, and then $\gamma$ is a continuous map. Since $\varphi^{-1}(A^*) = \pi(\gamma^{-1}(A^*)) = \pi(\overline{A})$ for each $A \in Z(\mathcal{LMC}(S))$, $\varphi$ is a continuous function.

$\Gamma$ is a Hausdorff space, because if $\mathcal{L}, \mathcal{K} \in \Gamma$ and $\mathcal{L} \neq \mathcal{K}$, then $\overline{\mathcal{L}} \bigcap \overline{\mathcal{K}} = \emptyset$. We can find open sets $V$ and $U$ in $S^{\mathcal{LMC}}$ such that $\overline{\mathcal{L}} \subseteq V \subseteq (\overline{\mathcal{K}})^c$, $\overline{\mathcal{K}} \subseteq U \subseteq (\overline{\mathcal{L}})^c$ and $U \cup V = S^{\mathcal{LMC}}$. Therefore there exist $A, B \in Z(\mathcal{LMC}(S))$ such that $\overline{\mathcal{K}} \subseteq U \subseteq \overline{A} \subseteq (\overline{\mathcal{L}})^c$ and $\overline{\mathcal{L}} \subseteq V \subseteq \overline{B} \subseteq (\overline{\mathcal{K}})^c$. Now let $\widetilde{p} \in S^{\mathcal{LMC}}$ then $\widetilde{p} \in \overline{A}$ or $\widetilde{p} \in \overline{B}$ and so $A^* \cup B^* = \Gamma$. This prove $(A^*)^c \cap (B^*)^c = \emptyset$. It is obvious that $\mathcal{L} \in (A^*)^c$ and $\mathcal{K} \in (B^*)^c$. $\varphi$ is continuous, one-to-one and onto therefore $\varphi$ is homeomorphism and $\Gamma \approx S^{\mathcal{LMC}}/R$. $\square$

**Lemma 3.6.** *Let $R$ be an equivalence relation on $S^{\mathcal{LMC}}$ such that $S^{\mathcal{LMC}}/R$ is Hausdorff. Let $\Gamma = \{\cap [\widetilde{p}]_R : \widetilde{p} \in S^{\mathcal{LMC}}\}$. Then with the $\tau_*$ quotient topology on $\Gamma$, $\Gamma \approx S^{\mathcal{LMC}}/R$.*

**Proof :** It is obvious that $\varphi : S^{\mathcal{LMC}}/R \to \Gamma$ defined by $\varphi([\widetilde{p}]) = \cap [\widetilde{p}]$ is continuous with $\tau_*$ quotient topology on $\Gamma$. Also it is obvious that $\varphi$ is an onto and one to one map.

$\Gamma$ is Hausdorff under $\tau^*$ quotient topology, because if $\mathcal{L}, \mathcal{K} \in \Gamma$ and $\mathcal{L} \neq \mathcal{K}$, then there exist $\widetilde{p}, \widetilde{q} \in S^{\mathcal{LMC}}$ such that $\overline{\mathcal{L}} = [\widetilde{p}] \neq [\widetilde{q}] = \overline{\mathcal{K}}$. There exist open subsets $U$ and $V$ of $S^{\mathcal{LMC}}$ and there exist $A, B \in Z(\mathcal{LMC}(S))$ such that

$$[\widetilde{p}] \subseteq U \subseteq \pi^{-1}(\pi(U)) \subseteq \overline{A},$$
$$[\widetilde{q}] \subseteq V \subseteq \pi^{-1}(\pi(V)) \subseteq \overline{B}$$

and $U \cup V = S^{\mathcal{LMC}}$. Hence $A_* \cup B_* = \Gamma$ and so $(A_*)^c \cap (B_*)^c = \emptyset$, where $\cap [\widetilde{p}] \subseteq (B_*)^c$ and $\cap [\widetilde{q}] \subseteq (A_*)^c$. Therefore $\Gamma$ is Hausdorff and so $\varphi$ is a homeomorphism. $\square$

**Corollary 3.7.** *$\tau^*$ and $\tau_*$ are equal on $\Gamma$.*

**Proof:** It is obvious. $\square$

Now by Corollary 3.7, we can say:

**Definition 3.8.** Let $\Gamma \subseteq Z(\mathcal{LMC}(S))$ be a set of pure $z-$filters. The quotient topology on $\Gamma$ is the topology generated by $\{(A^*)^c : A \in Z(\mathcal{LMC}(S))\}$ or generated by $\{(A_*)^c : A \in Z(\mathcal{LMC}(S))\}$.

**Theorem 3.9.** *Let $\Gamma$ be a set of pure $z-$filters on $\mathcal{LMC}(S)$ with the quotient topology. There is an equivalence relation $R$ on $S^{\mathcal{LMC}}$ such that $S^{\mathcal{LMC}}/R$ is Hausdorff, $\Gamma = \{\cap [p]_R : p \in S^{\mathcal{LMC}}\}$, and $\Gamma \approx S^{\mathcal{LMC}}/R$ if and only if both of the following statements hold:*

*(a) Given any topological choice function $f$ for $\Gamma$, there is a finite subfamily $\mathcal{F}$ of $\Gamma$ such that $\overline{S} = \bigcup_{\mathcal{L} \in \mathcal{F}} (\overline{f(\mathcal{L})})^\circ$,*

*(b) Given distinct $\mathcal{L}$ and $\mathcal{K}$ in $\Gamma$, there exist $A \in \mathcal{L}^\circ$ and $B \in \mathcal{K}^\circ$ such that whenever $\mathcal{C} \in \Gamma$, for each $H, T \in Z(\mathcal{LMC}(S))$ either $\overline{S} - (\overline{A})^\circ \subseteq (\overline{H})^\circ$ then $H \in \mathcal{C}$ or $\overline{S} - (\overline{B})^\circ \subseteq (\overline{T})^\circ$ then $T \in \mathcal{C}$.*

**Proof:** (Necessity) Statement ($a$) holds by Theorem 3.2. To establish statement ($b$) let $\mathcal{L}$ and $\mathcal{K}$ be distinct members of $\Gamma$. Note that $\Gamma$ is compact and Hausdorff, and $\gamma : S^{\mathcal{LMC}} \to \Gamma$ defined by $\gamma(\widetilde{p}) = \bigcap [\widetilde{p}]_R$ is a continuous map. Pick disjoint basic neighborhoods $(F^*)^c$ and $(G^*)^c$ of $\mathcal{L}$ and $\mathcal{K}$, respectively. Then $F^* \bigcup G^* = \Gamma$. Since $\gamma$ is continuous, we have $\overline{\mathcal{L}} \subseteq \gamma^{-1}((F^*)^c)$ and $\overline{\mathcal{K}} \subseteq \gamma^{-1}((G^*)^c)$. Therefore there exist $A \in \mathcal{L}$ and $B \in \mathcal{K}$ such that $\overline{\mathcal{L}} \subseteq (\overline{A})^\circ \subseteq \overline{A} \subseteq \gamma^{-1}((F^*)^c)$ and $\overline{\mathcal{K}} \subseteq (\overline{B})^\circ \subseteq \overline{B} \subseteq \gamma^{-1}((G^*)^c)$. Hence $\gamma^{-1}(F^*) \subseteq \overline{S} - (\overline{A})^\circ$ and $\gamma^{-1}(G^*) \subseteq \overline{S} - (\overline{B})^\circ$. This conclude that $\mathcal{C} \in \Gamma = F^* \bigcup G^*$ and so $\mathcal{C} \in F^*$ or $\mathcal{C} \in G^*$. If $\mathcal{C} \in F^*$ then $\overline{\mathcal{C}} \subseteq \gamma^{-1}(F^*) \subseteq \overline{S} - (\overline{A})^\circ$. Hence for each $H \in Z(\mathcal{LMC}(S))$, if $\overline{S} - (\overline{A}) \subseteq (\overline{H})^\circ$ then $H \in \mathcal{C}$.

(Sufficiency) Assume that ($a$) and ($b$) hold and observe that statement ($b$) of Theorem 3.2 follows. Suppose that

$R = \{(\widetilde{p}, \widetilde{q}) : \widetilde{p}, \widetilde{q} \in S^{\mathcal{LMC}}$ and that there is some $\mathcal{L} \in \Gamma$ such that $\mathcal{L} \subseteq \widetilde{p} \bigcap \widetilde{q}\}$.

we have as in Theorem 3.2 that $\Gamma = \{\bigcap [p]_R : p \in S^{\mathcal{LMC}}\}$. To complete the proof it suffices to show that R is closed in $S^{\mathcal{LMC}} \times S^{\mathcal{LMC}}$. (For then, since $S^{\mathcal{LMC}}$ is a compact Hausdorff space, $S^{\mathcal{LMC}}/R$ is Hausdorff and hence Lemma 3.5 applies.)

To this end let $(\widetilde{p}, \widetilde{q}) \in S^{\mathcal{LMC}} \times S^{\mathcal{LMC}} - R$. By Theorem 3.2 ($a'$), pick $\mathcal{L}$ and $\mathcal{K}$ in $\Gamma$ such that $\mathcal{L} \subseteq \widetilde{p}$ and $\mathcal{K} \subseteq \widetilde{q}$. Since $(\widetilde{p}, \widetilde{q}) \notin R$, $\mathcal{L} \neq \mathcal{K}$. Pick $A \in \mathcal{L}^\circ$ and $B \in \mathcal{K}^\circ$ such that for each $H, T \in Z(\mathcal{LMC}(S))$ either $\overline{S} - (\overline{A})^\circ \subseteq (\overline{H})^\circ$ then $H \in \mathcal{C}$ or $\overline{S} - (\overline{B})^\circ \subseteq (\overline{T})^\circ$ then $T \in \mathcal{C}$, where $\mathcal{C} \in \Gamma$. Thus $(\overline{A})^\circ \times (\overline{B})^\circ$ is a neighborhood of $(\widetilde{p}, \widetilde{q})$ in $S^{\mathcal{LMC}} \times S^{\mathcal{LMC}}$ which misses R. (If $(\mathcal{U}, \mathcal{V}) \in (\overline{A})^\circ \times (\overline{B})^\circ$ and $\mathcal{C} \in \Gamma$ with $\mathcal{C} \subseteq \mathcal{U} \cap \mathcal{V}$ and for each $H, T \in Z(\mathcal{LMC}(S))$ either $\overline{S} - (\overline{A})^\circ \subseteq (\overline{H})^\circ$ then $H \in \mathcal{C} \subseteq \mathcal{U}$, and so $\mathcal{U} \in (\overline{A})^\circ \cap (\bigcap_{(\overline{H})^\circ \supseteq \overline{S} - (\overline{A})^\circ} (\overline{H})) = (\overline{A})^\circ \cap (\overline{S} - (\overline{A})^\circ) = \emptyset$ is a contradiction.) □

**Definition 3.10.** $\Gamma$ is a quotient of $S^{\mathcal{LMC}}$ if and only if $\Gamma$ is a set of pure $z-$filters on $\mathcal{LMC}(S)$ with the quotient topology satisfying statements ($a$) and ($b$) of Theorem 3.9.

If $\Gamma$ is a quotient of $S^{\mathcal{LMC}}$, then by Theorem 3.2, for each $\widetilde{p} \in S^{\mathcal{LMC}}$ there is a unique $\mathcal{L} \in \Gamma$ such that $\mathcal{L} \subseteq \widetilde{p}$. Consequently, the function $\gamma$ below is well defined.

**Definition 3.11.** Let $\Gamma$ be a quotient of $S^{\mathcal{LMC}}$. Define $\gamma : S^{\mathcal{LMC}} \to \Gamma$ by $\gamma(\widetilde{p}) \subseteq \widetilde{p}$. Define $e : S \to \Gamma$ by $e(s) = \gamma(\widehat{s})$.

It is obvious that if $\{\widehat{s} : s \in S\} \subseteq \Gamma$ then $e(s) = \widehat{s}$.

**Corollary 3.12.** *Let $\Gamma$ be a quotient of $S^{\mathcal{LMC}}$. Then $\gamma$ is a quotient map.*

**Proof:** Let $R$, $\varphi : S^{\mathcal{LMC}}/R \to \Gamma$, and $\pi : S^{\mathcal{LMC}} \to S^{\mathcal{LMC}}/R$ be as in the proof of Lemma 3.5. By Theorem 3.9, $S^{\mathcal{LMC}}/R$ is Hausdorff. Thus $\varphi$ is a homeomorphism. Since $\pi$ is a quotient map and $\gamma = \varphi \circ \pi$, $\gamma$ is a quotient map. □

**Definition 3.13.** Let $\Gamma$ be a quotient of $S^{\mathcal{LMC}}$. If for each $\mathcal{A}$ and $\mathcal{B}$ in $\Gamma$ there is some $\mathcal{C} \in \Gamma$ with $\mathcal{C} \subseteq \mathcal{A} \odot \mathcal{B}$, define $\dotplus$ on $\Gamma$ by $\mathcal{A} \dotplus \mathcal{B} \in \Gamma$ and $\mathcal{A} \dotplus \mathcal{B} \subseteq \mathcal{A} \odot \mathcal{B}$.

**Lemma 3.14.** *Let $(\Gamma, \dot{+})$ be as in Definition 3.13. Then the following diagram commutes:*

$$\begin{array}{ccc} S^{\mathcal{LMC}} \times S^{\mathcal{LMC}} & \xrightarrow{\odot} & S^{\mathcal{LMC}} \\ \gamma \times \gamma \downarrow & & \downarrow \gamma \\ \Gamma \times \Gamma & \xrightarrow{\dot{+}} & \Gamma. \end{array}$$

**Proof:** We first observe that $\dot{+}$ is well defined. Indeed, if $\mathcal{C}$ and $\mathcal{D}$ are in $\Gamma$ and $\mathcal{C} \subseteq \mathcal{A} \odot \mathcal{B}$ and $\mathcal{D} \subseteq \mathcal{A} \odot \mathcal{B}$, then $\overline{\mathcal{A} \odot \mathcal{B}} \subseteq \overline{\mathcal{C}} \cap \overline{\mathcal{D}}$. But by condition $(b)$ of Theorem 3.9, if $\overline{\mathcal{C}} \cap \overline{\mathcal{D}} \neq \emptyset$, then $\mathcal{C} = \mathcal{D}$.

Now let $\widetilde{p}$ and $\widetilde{q}$ be in $S^{\mathcal{LMC}}$. Let $\mathcal{A} = \gamma(\widetilde{p})$, $\mathcal{B} = \gamma(\widetilde{q})$, $\mathcal{C} = \gamma(\widetilde{p} \odot \widetilde{q})$ and $\mathcal{D} = \gamma(\widetilde{p}) \dot{+} \gamma(\widetilde{q})$. Since $\widetilde{p} \odot \widetilde{q} \in \overline{\mathcal{A} * \mathcal{B}}$ so Lemma 2.18 implies that $\widetilde{p} \odot \widetilde{q} \in \overline{\mathcal{A} \odot \mathcal{B}}$. It is obvious that $\mathcal{D} \subseteq \mathcal{A} \odot \mathcal{B}$ so $\widetilde{p} \odot \widetilde{q} \in \overline{\mathcal{D}}$ also $\widetilde{p} \odot \widetilde{q} \in \overline{\mathcal{C}}$. Hence $\overline{\mathcal{C}} \cap \overline{\mathcal{D}} \neq \emptyset$ and this implies $\mathcal{C} = \mathcal{D}$. □

**Theorem 3.15.** *Let $(\Gamma, \dot{+})$ be as in Definition 3.13. Then*
*(1) $e$ is a continuous homomorphism from $S$ to $\Gamma$,*
*(2) $\Gamma$ is a right topological semigroup,*
*(3) $e[S]$ is dense in $\Gamma$, and*
*(4) for every $s \in S$, the function $\lambda_{e(s)}$ is continuous.*

**Proof:** (1) It is obvious that the following diagram commutes:

$$\begin{array}{ccc} S & \xrightarrow{\varepsilon} & S^{\mathcal{LMC}} \\ & e \searrow \quad \swarrow \gamma & \\ & \Gamma & \end{array}$$

(i.e. $\gamma \circ \varepsilon = e$). $\gamma$ and $\varepsilon$ are continuous so $e$ is continuous. By Lemma 3.14 $\gamma$ is a homomorphism, therefore $e$ is a continuous homomorphism.

(2) Let $\mathcal{L}$, $\mathcal{K}$ and $\mathcal{C}$ be in $\Gamma$. $\gamma$ is onto, therefore there are $\widetilde{p}$, $\widetilde{q}$ and $\widetilde{r}$ in $S^{\mathcal{LMC}}$ such that $\gamma(\widetilde{p}) = \mathcal{L}$, $\gamma(\widetilde{q}) = \mathcal{K}$ and $\gamma(\widetilde{r}) = \mathcal{C}$. By Lemma 3.14, we have

$$\begin{aligned} \mathcal{L} \dot{+} (\mathcal{K} \dot{+} \mathcal{C}) &= \gamma(\widetilde{p}) \dot{+} (\gamma(\widetilde{q}) \dot{+} \gamma(\widetilde{r})) \\ &= \gamma(\widetilde{p}) \dot{+} (\gamma(\widetilde{q} \odot \widetilde{r})) \\ &= \gamma(\widetilde{p} \odot (\widetilde{q} \odot \widetilde{r})) \\ &= \gamma((\widetilde{p} \odot \widetilde{q}) \odot \widetilde{r}) \\ &= \gamma(\widetilde{p} \odot \widetilde{q}) \dot{+} \gamma(\widetilde{r}) \\ &= (\gamma(\widetilde{p}) \dot{+} \gamma(\widetilde{q})) \dot{+} \gamma(\widetilde{r}) \\ &= (\mathcal{L} \dot{+} \mathcal{K}) \dot{+} \mathcal{C}. \end{aligned}$$

and so $\Gamma$ is a semigroup.

To see that $\dot{+}$ is right continuous, let $\mathcal{L} \in \Gamma$ and pick $\widetilde{p} \in \overline{\mathcal{L}} = \gamma^{-1}(\{\mathcal{L}\})$ then we have $r_\mathcal{L} \circ \gamma = \gamma \circ r_{\widetilde{p}}$. We know $r_{\widetilde{p}}$ and $\gamma$ are continuous. Thus $r_\mathcal{L}$ is continuous.

(3) By (1), we have $e = \gamma \circ \varepsilon$ and so $\Gamma = \gamma(\overline{\varepsilon(S)}) \subseteq \overline{\gamma(\varepsilon(S))} = \overline{\gamma \circ \varepsilon(S)} = \overline{e[S]}$.

(4) Let $s \in S$. It is obvious that $\lambda_{e(s)} \circ \gamma = \gamma \circ \lambda_{\varepsilon(s)}$. Since $\lambda_{\varepsilon(s)}$ and $\gamma$ are continuous so $\lambda_{e(s)}$ is continuous. $\square$

**Theorem 3.16.** *Let $\Gamma$ be a set of pure $z-$filters with the quotient topology. The following statements are equivalent:*

*(a) There exists a continuous function $h : S \longrightarrow \Gamma$ such that $\Gamma$ is a Hausdorff compact space, $h[S]$ is dense in $\Gamma$ and*

*(i) for each $s \in S$, $s \in \bigcap h(s)$ and*

*(ii) for distinct $\mathcal{L}$ and $\mathcal{K}$ in $\Gamma$, there exists $B \in \mathcal{K}^\circ$ such that for each $A \in Z(\mathcal{LMC}(S))$ if $\overline{S} - (\overline{B})^\circ \subseteq (\overline{A})^\circ$ then $A \in \mathcal{L}$.*

*(b) $\Gamma$ is quotient of $S^{\mathcal{LMC}}$.*

*(c) $\Gamma$ is quotient of $S^{\mathcal{LMC}}$, $\Gamma$ is a Hausdorff compact space and $e[S]$ is dense in $\Gamma$.*

**Proof:** To see that (a) implies (b), we must show that condition (a) and (b) of Theorem 3.9 hold. Let $f$ be a topological choice function for $\Gamma$. Then $\{\gamma((\overline{f(\mathcal{L})})^\circ) : \mathcal{L} \in \Gamma\}$ is an open cover of $\Gamma$ so pick finite $\mathcal{F} \subseteq \Gamma$ such that $\Gamma \subseteq \bigcup_{\mathcal{L} \in \mathcal{F}} \gamma((\overline{f(\mathcal{L})})^\circ)$. To see that $\overline{S} \subseteq \bigcup_{\mathcal{L} \in \mathcal{F}} \gamma((\overline{f(\mathcal{L})})^\circ)$, let $s \in S$. Pick $\mathcal{L} \in \mathcal{F}$ such that $h(s) \in \gamma((\overline{f(\mathcal{L})})^\circ)$. Then $f(\mathcal{L}) \in h(s)$, since $s \in \bigcap h(s)$ so $s \in f(\mathcal{L})$. Hence $S \subseteq \bigcup_{\mathcal{L} \in \mathcal{F}} (\overline{f(\mathcal{L})})^\circ$ and this implies that $\overline{S} \subseteq \bigcup_{\mathcal{L} \in \mathcal{F}} (\overline{f(\mathcal{L})})^\circ$.

The proof that condition (b) of Theorem 3.9 holds can be taken nearly verbatim from the proof of Theorem 3.9.

To see that (b) implies (c) observe that by Theorem 3.9, $\Gamma$ is a compact Hausdorff space. By Theorem 3.15, $e$ is continuous and $\overline{e[S]} = \Gamma$.

That (c) implies (a) is trivial. $\square$

Department of Mathematics, Faculty of science, Shahed University of Tehran P. O. Box: 18151-159, Tehran, Iran

*E-mail address*: `akbari@shahed.ac.ir`